\documentclass[a4paper,12pt]{article}

\usepackage{amsthm}
\usepackage{amsmath,amssymb,latexsym,amsfonts,mathrsfs}
\usepackage[dvips]{graphicx}

\theoremstyle{plain}
\newtheorem{Thm}{Theorem}[section]

\newtheorem{Lem}[Thm]{Lemma}

\newtheorem{Cor}[Thm]{Corollary}

\theoremstyle{definition}

\newcommand{\Proof}[2][Proof]{\begin{proof}[{#1}] #2 \end{proof}}

\renewcommand{\d}{{\rm d}} 

\newcommand{\tend}[2]{\mathrel{\mathop{\longrightarrow}\limits^{#1}_{#2}}}

\renewcommand{\bar}{\overline}
\newcommand{\un}[1]{\underline{#1}}

\numberwithin{equation}{section}

\makeatletter
\renewcommand\section{\@startsection {section}{1}{\z@}%
                                   {-3.5ex \@plus -1ex \@minus -.2ex}%
                                   {2.3ex \@plus.2ex}%
                                   {\normalfont\large\bf}}
\makeatother

\makeatletter
\renewcommand\subsection{\@startsection {subsection}{1}{\z@}%
                                   {-3.5ex \@plus -1ex \@minus -.2ex}%
                                   {2.3ex \@plus.2ex}%
                                   {\normalfont\normalsize\bf}}
\makeatother

\makeatletter
\@addtoreset{footnote}{page}
\makeatother

\newcommand{\rbra}[1]{\!\left( #1 \right)} 
\newcommand{\cbra}[1]{\!\left\{ #1 \right\}} 
\newcommand{\ceil}[1]{\!\left\lceil #1 \right\rceil} 

\newcommand{\bR}{\ensuremath{\mathbb{R}}}

\newcommand{\bZ}{\ensuremath{\mathbb{Z}}}

\newcommand{\cB}{\ensuremath{\mathcal{B}}}

\begin{document}

\begin{center}
{\Large \bf 
Entropy of random chaotic interval map with noise which causes coarse-graining 
}
\end{center}
\begin{center}
Kouji \textsc{Yano}\footnote{
Graduate School of Science, Kyoto University, Kyoto, Japan.}\footnote{
Research partially supported by KAKENHI (20740060), by KAKENHI (24540390) and by Inamori Foundation.}
\end{center}
\begin{center}
{\small \today}
\end{center}
\bigskip

\begin{abstract}
A random chaotic interval map with noise which causes coarse-graining 
induces a finite-state Markov chain. 
For a map topologically conjugate to a piecewise-linear map with the Lebesgue measure 
being ergodic, 
we prove that the Shannon entropy for the induced Markov chain 
possesses a finite limit as the noise level tends to zero. 
In most cases, the limit turns out 
to be strictly greater than the Lyapunov exponent of the original map without noise. 
\end{abstract}



\section{Introduction}

For the study of random mapping dynamics, Lyapunov exponents play key roles. 
Furstenberg--Kesten \cite{MR0121828} 
proved convergence of upper Lyapunov exponent for products of independent random matrices 
(see also Bougerol--Lacroix \cite{MR886674}). 
Diaconis--Freedman \cite{MR1669737} proved 
almost sure convergence of the backward iteration 
if the random mapping is contracting on the average. 
Steinsaltz \cite{MR1908881} proved 
almost sure convergence of the backward iteration 
for random logistic maps 
under the assumption that the averaged Lyapunov exponent is negative. 

Matsumoto--Tsuda \cite{MR711470} observed that 
the numerical KS entropy for a modified BZ map with noise 
may fall below that for the original map without noise, 
and called this phenomenon the {\em noise-induced order}. 
For mathematical results, 
Sumi \cite{MR2747724} proved that 
the chaos disappears 
for most of random complex dynamical systems for rational chaotic maps. 

In order to study 
how a (non-random) mapping dynamics is affected by a noise, 
it may be useful to study how the Lyapunov exponent is related 
to some entropies for random chaotic maps. 
Ara{\'u}jo--Tahzibi \cite{MR2134078} 
proved that the metric entropy of a random mapping dynamics, 
which was introduced by Kifer \cite{MR1015933} via its skew product realization, 
falls below the KS entropy of the noise zero limit of the random mapping dynamics. 
Kozlov--Treshchev \cite{MR2347306} 
and Piftankin--Treschev \cite{MR2679766} 
proved that 
the coarse-graining Gibbs entropy converges to the KS entropy 
in the noise zero limit. 

In this paper, we study the noise zero limit of the entropy of random chaotic maps 
through an approach which is different from all the above results. 

Let $ f $ be a chaotic map on the interval $ [0,1] $ 
with invariant probability measure $ \mu $ 
and consider a device which is designed to return $ f(x) $ as output 
if input is $ x $ and if there is no noise. 
Suppose there is a noise which affects the device in such a way as coarse-graining the states; 
more precisely, 
the states are clustered into the set of subintervals 
$ \Delta = \{ A^{(1)},\ldots,A^{(N)} \} $ 
equivolume with respect to $ \mu $, 
and, if input is $ n $ taken from $ \{ 1,\ldots,N \} $, 
the device picks a point $ U $ from the subinterval $ A^{(n)} $ at random 
with respect to $ \mu $ conditional on $ A^{(n)} $ 
and returns $ n' $ such that $ f(U) \in A^{(n')} $ as output. 
To iterate this procedure independently induces 
a Markov chain taking values in $ \{ 1,\ldots,N \} $. 

The purpose of this paper is to study 
the fine-graining limit as the noise level $ 1/N $ tends to zero 
of the Shannon entropy $ H_{\Delta}(f) $ for the induced Markov chain. 
We shall prove that 
$ \limsup H_{\Delta}(f) $ and $ \liminf H_{\Delta}(f) $ are invariants 
with respect to topological conjugate. 
We shall also prove that, 
for piecewise-linear map with the Lebesgue measure 
being ergodic, 
the fine-graining limit does exist and is obtained explicitly. 
It is remarkable that the limit is always no less, and, in most cases strictly greater, 
than the Lyapunov exponent $ \lambda(f) $ 
of the original (non-random) dynamical system $ (f,\mu) $. 

Let us give a small remark. 
Misiurewicz \cite{MR1002409} and \cite{MR1913288} 
studied continuity and discontinuity of topological entropies 
for piecewise monotone interval maps under perturbations 
preserving the number of pieces of monotonicity. 
He proved that the topological entropy for the skew tent maps is continuous. 
In a remarkable contrast, our fine-graining limit of the Shannon entropy 
for such a map is strictly greater than its Lyapunov exponent. 

We give another small remark. 
The induced Markov chain can always be realized as a random mapping dynamics. 
So one may want to adopt the Shannon entropy of the random mapping dynamics 
rather than that of the Markov chain. 
However, the former is not less than the latter, 
and, in addition, the way of such realizations is not unique; 
see Yano--Yasutomi \cite{MR2884814} and \cite{YanoYasuNova} for related results.

This paper is organized as follows. 
In Section \ref{sec: top inv}, we prepare notations 
of the finite-state Markov chain induced by coarse-graining. 
In Section \ref{sec: entropy}, 
we define $ H_{\Delta}(f) $ 
and prove that its fine-graining limits are invariants with respect to topological conjugate. 
Section \ref{sec: main} is devoted to the computation of the fine-graining limit. 
In Section \ref{sec: ex}, we examine the results in the case of skew tent maps.

\section{Random chaotic maps with noise which causes coarse-graining} \label{sec: top inv} 

Let $ f:[0,1] \to [0,1] $ be a measurable map 
with a unique absolutely continuous invariant probability measure $ \mu $ on $ [0,1] $ 
which is ergodic. 
For a positive integer $ N $, 
we call $ \Delta = \{ A^{(1)},\ldots,A^{(N)} \} $ 
an {\em equivolume partition} if 
$ \Delta $ consists of disjoint subintervals of $ [0,1] $ 
such that $ \bigcup_{n=1}^N A^{(n)} = [0,1] $ 
and $ \mu(A^{(n)}) = 1/N $ for $ n=1,\ldots,N $. 
We write $ \| \Delta \| = 1/N $, which will be called the {\em noise level}. 
Let $ U = (U^{(1)},\ldots,U^{(N)}) $ be a vector-valued random variable 
whose marginal $ U^{(n)} $ is distributed as $ \mu $ conditional on $ A^{(n)} $, 
i.e., 
\begin{align}
P \rbra{ U^{(n)} \in B } 
= \frac{\mu \rbra{ B \cap A^{(n)} }}{\mu \rbra{A^{(n)}} } 
= N \mu \rbra{ B \cap A^{(n)} } 
\quad \text{for} \ B \in \cB([0,1]) . 
\label{}
\end{align}
We do not require any assumption for the joint distribution 
among $ U^{(1)},\ldots,U^{(N)} $, 
because we only need the marginal distributions of $ U $. 
Let $ \pi^{\Delta}:[0,1] \to \{ 1,\ldots,N \} $ be the projection map such that 
\begin{align}
\pi^{\Delta}[x] = n 
\ \text{if and only if} \ 
x \in A^{(n)}. 
\label{}
\end{align}
We define a random map $ f^{\Delta} $ from $ [0,1] $ to itself by 
\begin{align}
f^{\Delta}(x) = f(U^{(\pi^{\Delta}[x])}) . 
\label{eq: fDelta}
\end{align}
We define a random map $ F^{\Delta} $ 
from $ \{ 1,\ldots,N \} $ to itself by 
\begin{align}
F^{\Delta}(n) = \pi^{\Delta}[f(U^{(n)})] . 
\label{eq: FDelta}
\end{align}
We note that 
\begin{align}
\pi^{\Delta}[f^{\Delta}(x)] = F^{\Delta}(\pi^{\Delta}[x]) . 
\label{}
\end{align}
For $ n,n' = 1,\ldots,N $, we write 
\begin{align}
p_{\Delta}(n'|n) 
:= P(F^{\Delta}(n)=n') 
= N \cdot \mu \rbra{ f^{-1} \rbra{A^{(n')}} \cap A^{(n)} } . 
\label{eq: pDelta}
\end{align}

We are now interested in the orbit of the iteration of the random maps repeated independently. 
Let $ (U_t)_{t=1,2,\ldots} $ be a sequence of independent copies of $ U $. 
Then we obtain the random maps $ (f_t^{\Delta})_{t=1,2,\ldots} $ 
and $ (F_t^{\Delta})_{t=1,2,\ldots} $ 
from \eqref{eq: fDelta} and \eqref{eq: FDelta}. 
Let $ x_0 $ be a random variable taking values in $ [0,1] $ 
and being independent of $ (U_t)_{t=1,2,\ldots} $ 
which obeys the law $ \mu $. 
Set $ X_0 = \pi^{\Delta}[x_0] $, 
which is thus distributed uniformly on $ \{ 1,\ldots,N \} $. 
We define $ (x_t)_{t=1,2,\ldots} $ and $ (X_t)_{t=1,2,\ldots} $ recursively by 
\begin{align}
x_t = f_t^{\Delta}(x_{t-1}) 
, \quad t=1,2,\ldots 
\label{}
\end{align}
and 
\begin{align}
X_t = F_t^{\Delta}(X_{t-1}) 
, \quad t=1,2,\ldots 
\label{}
\end{align}
We note that 
\begin{align}
X_t = \pi^{\Delta}[x_t] 
, \quad t=1,2,\ldots 
\label{}
\end{align}
and it is immediate that 
$ (X_t)_{t=1,2,\ldots} $ is a time-homogeneous Markov chain. 
Its transition probability is given as 
\begin{align}
P(X_t=n'|X_{t-1}=n) = p_{\Delta}(n'|n) 
\label{}
\end{align}
for $ n,n' = 1,\ldots,N $ and $ t=1,2,\ldots $, 
and its stationary distribution is the uniform distribution: 
\begin{align}
\mu_{\Delta}(n) = \frac{1}{N} 
, \quad n=1,\ldots,N 
\label{}
\end{align}

\section{Entropy of the induced Markov chain} \label{sec: entropy}

We denote the Shannon entropy of the induced Markov chain $ (X_t)_{t=1,2,\ldots} $ 
by 
\begin{align}
H_{\Delta}(f) = \sum_{n,n'=1}^N \mu_{\Delta}(n) \phi \rbra{ p_{\Delta}(n'|n) } , 
\label{}
\end{align}
where 
\begin{align}
\phi(t) = - t \log t 
\ (t>0) 
, \quad 
\phi(0)=0. 
\label{}
\end{align}
Now we write its fine-graining limits as the noise level $ \| \Delta \| = 1/N $ tends to zero by 
\begin{align}
\bar{H}(f) = \limsup_{\| \Delta \| \to 0} H_{\Delta}(f) 
, \quad 
\un{H}(f) = \liminf_{\| \Delta \| \to 0} H_{\Delta}(f) . 
\label{eq: def Hf}
\end{align}

\begin{Thm} \label{thm: inv}
Suppose that $ f:[0,1] \to [0,1] $ has 
a unique absolutely continuous invariant probability measure $ \mu $ on $ [0,1] $ 
which is ergodic. 
Then the fine-graining limits $ \bar{H}(f) $ and $ \un{H}(f) $ are 
invariants with respect to topological conjugate. 
\end{Thm}

\Proof{
Let $ C:[0,1] \to [0,1] $ be a homeomorphism 
and write $ g = C \circ f \circ C^{-1} $. 
Then the interval map $ g $ also has 
the unique absolutely continuous invariant probability measure 
given as $ \nu := \mu \circ C^{-1} $ 
which is ergodic. 
For any partition $ \Delta=\{ A^{(1)},\ldots,A^{(N)} \} $ equivolume with respect to $ \mu $, 
the partition $ C(\Delta) = \{ C(A^{(1)}),\ldots,C(A^{(N)}) \} $ 
is equivolume with respect to $ \nu $. 
Let us denote the transition probability $ p_{\Delta} $ 
for the dynamical system $ (f,\mu) $ and the equivolume partition $ \Delta $ 
as is defined in \eqref{eq: pDelta}, 
and write $ q_{C(\Delta)} $ for its counterpart for the dynamical system $ (g,\nu) $ 
and the equivolume partition $ C(\Delta) $. 
It is then obvious that 
\begin{align}
p_{\Delta}(n'|n) 
=& N \cdot \mu \rbra{ f^{-1} \rbra{A^{(n')}} \cap A^{(n)} } 
\label{} \\
=& N \cdot \nu \rbra{ g^{-1} \rbra{C(A^{(n')})} \cap C(A^{(n)}) } 
\label{} \\
=& q_{C(\Delta)}(n'|n) . 
\label{}
\end{align}
Now we obtain 
\begin{align}
H_{\Delta}(f) = H_{C(\Delta)}(g) . 
\label{}
\end{align}
If $ \Delta $ varies all the equivolume partitions, so does $ C(\Delta) $. 
Therefore, immediately from the definition \eqref{eq: def Hf}, we obtain 
\begin{align}
\bar{H}(f) = \bar{H}(g) 
, \quad 
\un{H}(f) = \un{H}(g) . 
\label{}
\end{align}
The proof is now complete. }

\section{Existence of fine-graining limits for piecewise-linear maps} \label{sec: main} 

Let $ f $ be an interval map 
with a unique absolutely continuous invariant probability measure $ \mu $ on $ [0,1] $ 
which is ergodic. 
Suppose that $ f $ is piecewise $ C^1 $, i.e., 
there exists a finite partition of $ [0,1] $, say $ 0=a_0<a_1<\cdots<a_{r-1}<a_r=1 $, 
such that 
the restriction of $ f $ on each subinterval $ [a_{i-1},a_i] $ 
can be extended to a $ C^1 $ map defined on an open interval including $ [a_{i-1},a_i] $. 
The Lyapunov exponent of $ f $ is defined as 
\begin{align}
\lambda(f) = \int_0^1 \log |f'(x)| \mu(\d x) . 
\label{}
\end{align}
Let us write 
\begin{align}
\cbra{x} = \min \{ x + n : n \in \bZ , \ x+n \ge 0 \} . 
\label{}
\end{align}

\begin{Thm} \label{thm: main} 
Suppose that $ \mu $ is the Lebesgue measure on $ [0,1] $. 
Suppose, in addition, that $ f $ is piecewise-linear, i.e., 
there exists a finite partition of $ [0,1] $, say $ 0=a_0<a_1<\cdots<a_{r-1}<a_r=1 $, 
such that $ f $ is linear on each subinterval $ E_i = (a_{i-1},a_i) $. 
Then one has 
\begin{align}
H(f) := \bar{H}(f) = \un{H}(f) = \lambda(f) + D(f) , 
\label{}
\end{align}
where $ D(f) $ is given as 
\begin{align}
D(f) = 2 \int_0^1 \frac{\rho(|f'(x)|)}{|f'(x)|} \d x 
\label{}
\end{align}
and the function $ \rho $ is defined as 
\begin{align}
\rho(m) =& 0 
\quad \text{if $ m \in \bZ $}, 
\label{} \\
\rho(m) =& \frac{1}{p} \sum_{n=1}^{p-1} \phi \rbra{ \frac{n}{p} } 
\quad \text{if $ m = \frac{q}{p} $: irreducible, $ p,q \in \bZ $ and $ p \ge 2 $}, 
\label{} \\
\rho(m) =& \frac{1}{4} 
\quad \text{if $ m $ is irrational}. 
\label{}
\end{align}
\end{Thm}

Combining Theorem \ref{thm: main} with Theorem \ref{thm: inv}, we obtain the following. 

\begin{Cor}
Suppose that $ f $ is a piecewise-$ C^1 $ map which is topologically conjugate 
to a piecewise-linear map $ g $ with the Lebesgue measure 
being the unique absolutely continuous invariant probability measure $ \mu $ on $ [0,1] $ 
which is ergodic. 
Then one has 
\begin{align}
\bar{H}(f) = \un{H}(f) \ge \lambda(f) . 
\label{eq: ineq}
\end{align}
Unless $ g' $ is integer valued, 
the inequality in \eqref{eq: ineq} is strict. 
\end{Cor}

Before proving Theorem \ref{thm: main}, we need the following lemma. 

\begin{Lem} \label{lem: expanding}
Suppose that $ \mu $ is the Lebesgue measure. Then the map $ f $ satisfies 
\begin{align}
|f'| \ge 1 
\quad \text{a.e.} 
\label{eq: expanding}
\end{align}
\end{Lem}

\Proof[Proof of Lemma \ref{lem: expanding}]{
Recall that the operator $ L:L^1([0,1]) \to L^1([0,1]) $ defined as 
\begin{align}
(L \varphi)(x) = \sum_{y:f(y)=x} \frac{1}{|f'(y)|} \varphi(y) 
\quad \text{for $ \varphi \in L^1([0,1]) $}. 
\label{}
\end{align}
is the Perron--Frobenius operator for the dynamical system $ (f,\mu) $, i.e., 
\begin{align}
\int_0^1 (L \varphi)(x) \psi(x) \d x 
= \int_0^1 \varphi(x) \psi(f(x)) \d x 
\label{}
\end{align}
holds for all $ \varphi \in L^1([0,1]) $ and all $ \psi \in L^{\infty }([0,1]) $. 
If we take $ \varphi(x) \equiv 1 $, we have, 
since $ \mu \circ f^{-1} = \mu $, 
\begin{align}
\int_0^1 (L 1)(x) \psi(x) \d x 
= \int_0^1 \psi(f(x)) \d x = \int_0^1 \psi(x) \d x , 
\label{}
\end{align}
and thus we obtain 
\begin{align}
(L 1)(x) = \sum_{y:f(y)=x} \frac{1}{|f'(y)|} = 1 
\quad \text{a.e.} 
\label{}
\end{align}
From this we obtain \eqref{eq: expanding}. 
}

Now we prove Theorem \ref{thm: main}. 

\Proof[Proof of Theorem \ref{thm: main}]{
Let $ \Delta = \{ A^{(1)},\ldots,A^{(N)} \} $ be an equivolume partition. 
Since $ \mu $ is the Lebesgue measure, we may assume that 
$ A^{(n)}=[x_{n-1},x_n) $ for $ n=1,\ldots,N-1 $ and $ A^{(N)}=[x_{N-1},x_N] $ 
where $ x_n = n/N $ for $ n=1,\ldots,N $. 
For $ i=1,\ldots,r $, 
let $ m_i = |f'(x)| $ for $ x \in E_i $. 
Since $ m_i \ge 1 $ by Lemma \ref{lem: expanding} 
and since $ f' $ is constant on $ E_i $,  
we may suppose without loss of generality that $ m = f'(x) \ge 1 $ for $ x \in E_i $, 
and, consequently, $ f $ is increasing on $ E_i $. 

Let $ i=1,\ldots,r $ and $ n=1,\ldots,N $ be fixed such that $ \bar{A^{(n)}} \subset E_i $. 
Let us write $ m $ simply for $ m_i $. 
We now have 
\begin{align}
\mu(f(A^{(n)})) = f(x_n) - f(x_{n-1}) 
= m (x_n-x_{n-1}) 
= \frac{m}{N} . 
\label{}
\end{align}
Let $ u $ and $ v $ be such that 
\begin{align}
f(x_{n-1}) \in A^{(u)} 
\quad \text{and} \quad 
f(x_n) \in A^{(v)} . 
\label{}
\end{align}
We then have 
\begin{align}
u = \ceil{ N f(x_{n-1}) } 
\quad \text{and} \quad 
v = \ceil{ N f(x_n) } 
\label{}
\end{align}
where 
\begin{align}
\ceil{x} = \min \{ n \in \bZ : n \ge x \} 
\quad \text{for $ x \in \bR $}. 
\label{}
\end{align}
Set $ B^{(i)} = A^{(i)} $ for $ i=u+1,\ldots,v-1 $ and set 
\begin{align}
B^{(u)} = [f(x_{n-1}),x_u) 
\quad \text{and} \quad 
B^{(v)} = [x_{v-1},f(x_n)] . 
\label{}
\end{align}
We then have $ f(A^{(n)}) = B^{(u)} \cup \cdots \cup B^{(v)} $, and hence 
\begin{align}
p_{\Delta}(n'|n) = 0 
\quad \text{for $ n' < u $ or $ n' > v $}. 
\label{}
\end{align}
Since $ \mu $ is the Lebesgue measure 
and since $ f $ is linear on $ A^{(n)} $, we see that 
\begin{align}
p_{\Delta}(n'|n) 
= \frac{\mu(B^{(n')})}{\mu(f(A^{(n)}))} 
\quad \text{for $ n'=u,\ldots,v $}. 
\label{}
\end{align}
Hence we obtain 
\begin{align}
p_{\Delta}(n'|n) 
=& 
\begin{cases}
1/m & \text{if $ n'=u+1,\ldots,v-1 $}, \\
a/m & \text{if $ n'=u $}, \\
b/m & \text{if $ n'=v $}, \\
0 & \text{otherwise}, 
\end{cases}
\label{}
\end{align}
where 
\begin{align}
a =& \ceil{ N f(x_{n-1}) } - N f(x_{n-1}) 
= 1 - \cbra{ N f(x_{n-1}) } 
, \label{} \\
b =& N f(x_n) - \ceil{ N f(x_n) } + 1 
= \cbra{ N f(x_n) }_+ , 
\label{}
\end{align}
where $ \cbra{x}_+ = \min \{ x+n : n \in \bZ , \ x+n>0 \} $. 
Noting that $ \phi(xy) = x \phi(y) + y \phi(x) $ for $ x,y \ge 0 $, we have 
\begin{align}
\sum_{n'} \phi \rbra{ p_{\Delta}(n'|n) } 
=& (v-u-1) \phi \rbra{ \frac{1}{m} } 
+ \phi \rbra{ \frac{a}{m} } + \phi \rbra{ \frac{b}{m} } 
\label{} \\
=& (v-u-1+a+b) \phi \rbra{ \frac{1}{m} } 
+ \frac{1}{m} \phi(a) + \frac{1}{m} \phi(b) 
\label{} \\
=& m \phi \rbra{ \frac{1}{m} } 
+ \frac{1}{m} \phi(a) + \frac{1}{m} \phi(b) 
\label{} \\
=& \log m + \frac{1}{m} \phi(a) + \frac{1}{m} \phi(b) . 
\label{eq: formula}
\end{align}

Let $ i=1,\ldots,r $ be fixed and return to write $ m_i $ instead of $ m $. 
We then have 
\begin{align}
\sum_{n:\bar{A^{(n)}} \in E_i} \phi(b) 
= \sum_{n:\bar{A^{(n)}} \in E_i} \phi(\cbra{N f(n/N)}_+) 
= \sum_{n:\bar{A^{(n)}} \in E_i} \phi(\cbra{m_in}) , 
\label{}
\end{align}
where we note that $ \phi(\cbra{x}_+) = \phi(\cbra{x}) $ 
because $ \phi(0)=\phi(1)=0 $. 
Let $ c_i(N) $ denote the number of $ n $'s such that $ n:\bar{A^{(n)}} \in E_i $. 
Then we see that 
\begin{align}
\frac{1}{c_i(N)} \sum_{n:\bar{A^{(n)}} \in E_i} \phi(b) 
\tend{}{N \to \infty } \rho(m_i) 
\label{eq: average}
\end{align}
by the following arguments: 
\begin{enumerate}
\item 
If $ m_i \in \bZ $, then $ \cbra{m_in} = 0 $ for all $ n $ 
so that we obtain \eqref{eq: average}. 
\item 
If $ m_i = q/p $: irreducible, $ p,q \in \bZ $ and $ p \ge 2 $, 
then the set $ \{ \cbra{m_in} : n=kp+1,kp+2,\ldots,(k+1)p \} $ 
coincides with the set $ \{ 0,1/p,2/p,\ldots,(p-1)/p \} $ 
for all $ k=0,1,\ldots $, 
so that we obtain \eqref{eq: average}. 
\item 
If $ m_i $ is irrational, then by Weyl's equidistribution theorem we obtain 
\begin{align}
\frac{1}{c_i(N)} \sum_{n:\bar{A^{(n)}} \in E_i} \phi(\cbra{m_in}) 
\tend{}{N \to \infty } 
\int_0^1 \phi(t) \d t = \frac{1}{4} , 
\label{}
\end{align}
so that we obtain \eqref{eq: average}. 
\end{enumerate}
In the same way, we obtain 
\begin{align}
\frac{1}{c_i(N)} \sum_{n:\bar{A^{(n)}} \in E_i} \phi(a) 
\tend{}{N \to \infty } \rho(m_i) . 
\label{}
\end{align}
Thus, using \eqref{eq: formula}, we obtain 
\begin{align}
& \sum_{i=1}^r \sum_{n:\bar{A^{(n)}} \subset E_i } 
\mu_{\Delta}(n) \sum_{n'} \phi \rbra{ p_{\Delta}(n'|n) } 
\label{} \\
=& \sum_{i=1}^r \frac{c_i(N)}{N} 
\cdot \frac{1}{c_i(N)} \sum_{n:\bar{A^{(n)}} \subset E_i } \phi \rbra{ p_{\Delta}(n'|n) } 
\label{} \\
\tend{}{N \to \infty }& 
\sum_{i=1}^r (a_i-a_{i-1}) \cdot \cbra{ \log m_i + \frac{\rho(m_i)}{m_i} +\frac{\rho(m_i)}{m_i} } 
\label{} \\
=& \int_0^1 \cbra{ \log |f'(x)| + 2 \frac{\rho(|f'(x)|)}{|f'(x)|} } \d x . 
\label{}
\end{align}

Note that if $ \bar{A^{(n)}} $ does not included in any $ E_i $, 
then $ \bar{A^{(n)}} $ contains at least one of the points $ a_1,\ldots,a_r $, 
so that the number of such $ n $'s is not greater than $ r $. 
Since 
\begin{align}
\sum_{n'} \phi \rbra{ p_{\Delta}(n'|n) } 
\le \sum_{n'} \phi \rbra{ \frac{1}{N} } = \log N , 
\label{}
\end{align}
we obtain 
\begin{align}
\sum_{i=1}^r \sum_{n:\bar{A^{(n)}} \not\subset E_i } 
\mu_{\Delta}(n) \sum_{n'} \phi \rbra{ p_{\Delta}(n'|n) } 
\le r \cdot \frac{1}{N} \cdot \log N 
\tend{}{N \to \infty } 0 . 
\label{}
\end{align}
Therefore, we conclude that 
\begin{align}
H_{\Delta}(f) 
= \sum_{i=1}^r \sum_{n,n'} \mu_{\Delta}(n) \phi \rbra{ p_{\Delta}(n'|n) } 
\tend{}{N \to \infty } 
\int_0^1 \cbra{ \log |f'(x)| + 2 \frac{\rho(|f'(x)|)}{|f'(x)|} } \d x , 
\label{}
\end{align}
which completes the proof. 
}

\section{Examples: skew tent maps} \label{sec: ex} 

For the illustration of Theorem \ref{thm: main}, 
we compute the difference $ D(f) $ 
for a {\em skew tent map} $ f:[0,1] \to [0,1] $, which is defined as 
\begin{align}
f(x) = 
\begin{cases}
mx & \text{if $ 0 \le x \le 1/m $}, \\
l(1-x) & \text{if $ 1/m < x \le 1 $} 
\end{cases}
\label{}
\end{align}
for some $ m>1 $ and $ l>1 $ such that $ 1/m+1/l=1 $. 
Note that the Lebesgue measure is the unique absolutely continuous invariant probability measure 
for $ f $ and is ergodic; see, e.g., Jetschke--Stiewe \cite{MR914048}. 

(1) Suppose that $ m=2 $. In this case, we have $ l=2 $ and hence we have $ D(f)=0 $. 
Note that this map $ f $ is topologically conjugate to the {\em logistic map} 
\begin{align}
g(x) = \frac{1}{4} x(1-x) 
\quad \text{for $ x \in [0,1] $} , 
\label{}
\end{align}
so that we have 
\begin{align}
H(g) = \lambda(g) = H(f) = \lambda(f) = \log 2. 
\label{}
\end{align}

(2) Suppose that $ m $ is rational and $ m \neq 2 $. 
We represent $ m = (p+q)/p $ as an irreducible fraction with $ p,q \in \bZ $, $ p,q \ge 1 $. 
In this case, we have $ l = (p+q)/q $ and hence we obtain 
\begin{align}
D(f) 
=& 2 \cdot \frac{1}{m} \cdot \frac{\rho(m)}{m} 
+ 2 \cdot \frac{1}{l} \cdot \frac{\rho(l)}{l} 
\label{} \\
=& 2 \frac{p}{(p+q)^2} \sum_{n=1}^{p-1} \phi \rbra{ \frac{n}{p} } 
+ 2 \frac{q}{(p+q)^2} \sum_{n=1}^{q-1} \phi \rbra{ \frac{n}{q} } , 
\label{}
\end{align}
where the summation $ \sum_{n=1}^{p-1} $ (resp. $ \sum_{n=1}^{q-1} $) 
is discarded if $ p=1 $ (resp. $ q=1 $). 

(3) Suppose that $ m $ is irrational. 
In this case, we see that $ l $ is also irrational and hence we obtain 
\begin{align}
D(f) 
=& 2 \cdot \frac{1}{m} \cdot \frac{1/4}{m} + 2 \cdot \frac{1}{l} \cdot \frac{1/4}{l} 
\label{} \\
=& \frac{m^2-2m+2}{2m^2} . 
\label{}
\end{align}

\def\cprime{$'$}

\end{document}